\documentclass[utf8,11pt,english]{aclart}

\usepackage{hyperref}

\numberwithin{subsection}{section}
\numberwithin{subsubsection}{section}
\numberwithin{paragraph}{section}
\numberwithin{equation}{section}

\makeatletter
\usepackage{natbib}

\def\bibliography#1{%
\if@filesw
\immediate\write\@auxout{\string\bibdata{#1}}%
\fi
\expandafter\input{\bbl@main@language bst.tex}%
\@input@{\jobname.bbl}
}
\makeatother

\usepackage{tikz}
\usetikzlibrary{matrix,arrows}

\title{Diophantine Geometry and Analytic Spaces}
\author{Antoine Chambert-Loir}
\address{Département de mathématiques d'Orsay \\ Bâtiment 425 \\
Faculté des sciences d'Orsay \\ Université Paris-Sud \\ F-91405 Orsay Cedex}
\email{Antoine.Chambert-Loir@math.u-psud.fr}

\def\A{\mathbf A}
\def\Q{\mathbf Q}
\def\R{\mathbf R}
\def\P{\mathbf P}
\def\C{\mathbf C}

\def\Z{\mathbf Z}

\let\ra\rightarrow

\let\hra\hookrightarrow

\let\eps\varepsilon \let\phi\varphi
\let\emptyset\varnothing
\def\abs#1{\left|#1\right|}
\def\norm#1{\left\|#1\right\|}

\def\an{{\text{\upshape an}}}
\def\div{\operatorname{div}}
\def\Card{\operatorname{Card}}
\def\Gal{\operatorname{Gal}}
\def\Pic{\operatorname{Pic}}
\def\Spec{\operatorname{Spec}}
\def\hPic{\mathop{\widehat{\operatorname{Pic}}\vphantom{Pic}}}
\def\hdeg{\mathop{\widehat{\operatorname{deg}}\vphantom{deg}}}

\def\gm{\mathbf G_{\mathrm m}}
\def\gmv{\mathbf G^\an_{\mathrm m,v}}
\def\gmvg{\mathbf G^{g,\an}_{\mathrm m,v}}

\begin{document}
\maketitle

\section{Introduction}

Diophantine Geometry can be roughly defined as the \emph{geometric}
study of \emph{diophantine equations}. Historically, and for most
mathematicians, those equations are polynomial equations
with integer coefficients and one seeks for integer, or rational,
solutions; generalizations to number fields come naturally.
However, it has been discovered in the \textsc{xix}\textsuperscript{th} century
that number fields  share striking
similarities with finite extensions of the field $k(t)$
of rational functions with coefficients in a field~$k$,
the analogy being the best when $k$ is a finite field.
From this point of view, rings of integers of number fields
are analogues of rings of regular functions on a regular curve.
Namely, both rings are Dedekind
(\emph{i.e.}, integrally closed, one-dimensional, Noetherian) domains.

When one studies Diophantine Geometry over number fields,
the geometric shape defined by the polynomial equations
over the complex numbers plays an obvious important role.
Be it sufficient to recall the statement of Mordell conjecture (proved
by~\cite{faltings83}): a Diophantine equation whose
associated complex shape is a compact Riemann surface of genus at least~$2$
has only finitely many solutions.
Over function fields, such a role can only be played by analytic
geometry over non-archimedean fields, a much more recent
theory than its complex counterpart.

The lecture I gave at the Bellairs Workshop in  Number Theory
was devoted to a survey of recent works in Diophantine
Geometry over function fields, where analytic geometry
over non-archimedean fields in the sense of~\cite{berkovich1990}
took a significant place.
Since this topic was not the main one of the conference, the talk
had been deliberately informal and
the present notes aim at maintaining this character, in the hope
that they will be useful for geometers of all obediences,
be it Diophantine, tropical, complex, non-archimedean...

I thank Matt Baker, Xander Faber, and Joe Rabinoff for
their comments during the talk, and for having sent me their
notes. I especially thank Xander Faber for his efficient cheerleading,
without which this survey would probably never have appeared.
I also thank the two referees for their many valuable suggestions.

At the time of the conference, 
I was affiliated to 
Université de Rennes~1
and a member of Institut universitaire de France.

\section{The standard height function}

In all the sequel, we fix a field~$F$ which can be, either
the field~$\Q$ of rational numbers (arithmetic case),
or the field~$k(T)$ of rational functions with coefficients in a given field~$k$
(geometric case). This terminology will be explained later.
We let $\bar F$ be an algebraic closure of~$F$.

The standard height function  $h\colon \P^n(\bar F)\ra\R$
is a function measuring the ``complexity'' of a point in projective space
with homogeneous coordinates in~$\bar F$.

We begin by describing it on the subset~$\P^n(F)$
of $F$-rational points.

\subsection{Arithmetic case}
Let $\mathbf x$ be a point in~$\P^n(\Q)$.
We may assume that its homogeneous coordinates $[x_0:\dots:x_n]$ are chosen
so as to be coprime integers; they are then well defined up
to a common sign and one defines
\begin{equation}
h(\mathbf x) = \log \max(\abs {x_0},\dots,\abs{x_n}).
\end{equation}

One first observes a \emph{finiteness property:} for any real number~$B$,
there are only finitely many points~$\mathbf x\in\P^n(\Q)$
such that $h(\mathbf x)\leq B$. Indeed, this bound gives only
finitely many possibilities for each coordinate.

The height function behaves well under morphisms. Let $f\colon\P^n\dashrightarrow \P^m$ be a rational map given by homogeneous forms  $f_0,\dots,f_m\in\Q[X_0,\dots,X_n]$ of degree~$d$, without common factor.
Its exceptional locus~$E$ is the closed subspace in~$\P^n$
defined by the simultaneous vanishing of all~$f_i$s. Then, one proves
easily that  there exists a constant~$c$ such that
$h(f(\mathbf x))\leq d h(\mathbf x)+c$,
for any point $\mathbf x\in\P^n(\Q)$ such that $\mathbf x\notin E$.
The converse inequality is more subtle and relies on the \emph{Nullstellensatz:}
\emph{Let $X$ be a closed subscheme of~$\P^n$
such that $X\cap E=\emptyset$; then, there exists a constant~$c_X$
such that $h(f(\mathbf x))\geq dh(\mathbf x)-c_X$
for any $\mathbf x\in X(\Q)$.}

To add on the subtlety behind these apparently simple estimates,
let me remark that it is easy, given an explicit map~$f$,
to write down an explicit acceptable constant~$c$;
however, giving an explicit constant~$c_X$ requires
a quite non-trivial statement called the
\emph{effective arithmetic Nullstellensatz};
see, for example, \citet*{dandrea-krick-sombra2011}
for a recent and sharp version.

\subsection{Geometric case}
Let now $\mathbf x$ be a point in~$\P^n(k(T))$.
Again, we may choose a system of homogeneous coordinates
$[x_0:\dots:x_n]$ of~$\mathbf x$ consisting of polynomials
in~$k[T]$ without common factors. Such a system is unique
up to multiplication by a common nonzero constant.
Let us define the height of~$\mathbf x$ by the formula
\begin{equation}
h(\mathbf x) = \max(\deg x_0,\dots,\deg x_n). \end{equation}

If the base field~$k$ is finite, then the height satisfies
a similar finiteness property as in the arithmetic case:
since there are only finitely many polynomials $f_i\in k[T]$
of given degree, the set of points~$\mathbf x\in \P^n(k(T))$
such that $h(\mathbf x)\leq B$ is finite, for any~$B$.

The height function has exactly the same properties with respect to morphism
as in the arithmetic case.

\subsection{Geometric interpretation (geometric case)}
In the geometric case, the height
can be given a \emph{geometric interpretation}, free
of homogeneous coordinates. Indeed, let $C$ be the projective line
over~$k$, that is, the unique projective regular
$k$-curve with function field~$F=k(T)$.
Any point $\mathbf x\in\P^n(F)$ can be interpreted
as a morphism $\phi_{\mathbf x}\colon C\ra\P^n$ of $k$-schemes.
When $\phi_{\mathbf x}$ is generically one-to-one, then
$h(\mathbf x)$ can be computed as the degree of the
the rational curve~$C$, as  embedded in~$\P^n$ through~$\phi_{\mathbf x}$.
In the general case, one has
\begin{equation}
h(\mathbf x)=
\deg \phi_{\mathbf x}^*\mathscr O(1), \end{equation}
that is, $h(\mathbf x)$ is the degree
of the pull-back to~$C$ of the tautological line bundle on~$\mathbf P^n$.

\subsection{Extension to the algebraic closure (geometric case)}
The previous geometric interpretation suggests a way to
define the height on the whole of $\mathbf P^n(\bar F)$.
Namely, let $E$ be a finite extension of~$F$;
it is the field of rational functions on a projective regular curve~$C_E$
defined over a finite extension of~$k$.
Any point $\mathbf x\in \P^n(E)$ can be interpreted as a morphism
$\phi_{\mathbf x}\colon C_E\ra \mathbf P^n$ and one sets
\begin{equation}
h(\mathbf x) = \frac1{[E:F]} \deg \phi_{\mathbf x}^*\mathscr O(1). \end{equation}
One checks that the right-hand-side of
this formula does not depend on the actual choice
of a finite extension~$E$ such that $\mathbf x\in\P^n(E)$,
thus defining a function $h\colon\P^n(\bar F)\ra\R$.

\subsection{Absolute values}
Using absolute values,
one can give a general definition of the standard height function, valid
for any finite extension of~$F$.

Recall that an absolute value on a field~$F$ is
a map $\abs\cdot\colon F\ra\R_{\geq 0}$ subject to
the following axioms:
$\abs 0=0$, $\abs 1=1$, $\abs {ab}=\abs a\abs b$
and $\abs{a+b}\leq \abs{a}+\abs b$
for any $a,b\in F$.
Two absolute values $\abs\cdot$ and $\abs\cdot'$
are said to be equivalent if there exists a positive
real number~$\lambda$ such that $\abs a'=\abs a^\lambda$
for any $a\in F$. The trivial absolute value on~$F$
is defined by $\abs a_0=1$ for any $a\in F^*$.

Let $M_F$ be the set of non-trivial absolute values of~$F$, up to equivalence.
Any class $v\in F$ possesses  a preferred, normalized, representative,
denoted $\abs\cdot_v$, so that the \emph{product formula} holds:
\begin{equation}
\prod_{v\in M_F} \abs a_v = 1 \qquad \text{for any $a\in F^\times$.}
\end{equation}
It connects the non-trivial absolute values (on the left-hand-side)
and the trivial one (on the right-hand-side).

The field $F=\Q$ possesses the usual archimedean absolute value,
denoted $\abs\cdot_\infty$.  Absolute values nonequivalent to that
satisfy the \emph{ultrametric inequality} $\abs{a+b}\leq\max(\abs a,\abs b)$,
and each of them is associated  to a prime number~$p$.
The corresponding normalized $p$-adic absolute value is characterized
by the equalities $\abs p_p=1/p$ and $\abs a_p=1$ for any integer~$a$
which is prime to~$p$.  Therefore, $M_\Q=\{\infty,2,3,5,7,\dots\}$.
A similar description applies to number fields, the normalized
ultrametric absolute values are in correspondence with the maximal
ideals of the ring of integers, while the archimedean absolute
values correspond to real or pair of conjugate complex embeddings
of the field.

All absolute values of the field $F=k(T)$  which are trivial
on~$k$ are ultrametric.
They correspond to the closed points of the projective line~$\P^1_k$
(whose field of rational functions is precisely~$F$).
For example, there is a unique absolute value on~$k(T)$ which maps
a polynomial~$P$ to $e^{\deg(P)}$, and it corresponds
to the point at infinity of $\P^1_k$.
Similarly, let $A\in k[T]$ be a monic, irreducible, polynomial;
there is a unique absolute value on~$k(T)$ which 
maps~$A$ to $e^{-\deg(A)}$ and maps to~$1$ any polynomial
which indivisible by~$A$; this absolute value corresponds
to the closed point~$(A)$ of the affine line.

More generally, if $E$ is a finite extension of~$F$,
the set~$M_E$ is naturally in bijection with the set
of closed points of the unique projective normal curve~$C_E$
with function field~$E$.

The product formula is nothing but the formula that claims
that the number of zeroes of a rational function on a curve
is equal to the number of poles (in both cases,
counted with multiplicity).

In this language, the height function on $\P^n(\bar F)$ can be defined as
\begin{equation}
h(\mathbf x) = \frac1{[E:F]}
\sum_{v\in M_E} \log \max (\abs{x_0}_v,\dots,\abs{x_n}_v)
,
\end{equation}
where $E$ is a finite extension of~$F$ and
$\mathbf x=[x_0:\dots:x_n]\in \P^n(E)$.

\subsection{Properties}
In the arithmetic case, or, in the geometric case
over a \emph{finite} base field~$k$, the height function satisfies
an important \emph{finiteness principle}, due to~\cite{northcott1950}:
for any real number~$B$ and any positive integer~$d$,
the set of points $\mathbf x\in\P^n(\bar F)$
such that $[F(\mathbf x):F]\leq d $ and $h(\mathbf x)\leq B$
is finite. Obviously, this property does not hold in
the geometric case, when the base field is infinite.

In all cases, the height function has a similar behavior
with respect to morphisms. \emph{Let $f\colon \P^n\dashrightarrow\P^m$
be a rational map defined by homogeneous polynomials $(f_0,\dots,f_m)$
of degree~$d$, without common factor; let $E\subset\P^n$ be the locus
defined by $f_0,\dots,f_m$. Then, there exists a constant~$c_f$
such that $h(f(\mathbf x))\leq d h(\mathbf x)+c_f$ for any $\mathbf x\in\P^n(\bar F)$ such that $\mathbf x\notin E$.
Let $X$ be a closed subscheme of~$\P^n$ such that $X\cap E=\emptyset$;
then, there exists  a real number~$c_X$ such that
$h(f(\mathbf x))\geq d h(\mathbf x)-c_X$ for any $\mathbf x\in X(\bar F)$.}

\section{Heights for line bundles, canonical heights}

\subsection{Heights for line bundles}
For applications, it is important  to understand precisely
the behavior of heights under morphisms.
This is embodied in the following fact, called the \emph{height machine}.
Let $\mathscr F(X(\bar F);\R)$ be the vector space
of real valued functions on~$X(\bar F)$, and let $\mathscr F_{\mathrm b}(X(\bar F);\R)$ be its subspace of bounded functions.
There exists a unique additive map
\[ h \colon \Pic(X)\otimes_\Z\R \ra \mathscr F(X(\bar F);\R)/\mathscr F_{\mathrm b}(X(\bar F);\R), \quad L\mapsto h_L \]
such that  for any closed embedding $f\colon X\hra \P^n_F$
of~$X$ into a projective space,
\[ h_{f^*\mathscr O(1)}  \equiv h_{\mathbf P^n}\circ f \pmod {\mathscr F_{\mathrm b}(X(\bar F);\R)}, \]
where we have denoted $h_{\mathbf P^n}$ the height on projective
space that we had define in the previous section.

Uniqueness comes from the fact that $\Pic(X)\otimes_\Z\R$ is generated
by line bundles of the form $f^*\mathscr O(1)$, for some
closed embedding~$f$. The existence follows from basic
properties of the height on projective spaces, namely
its behavior under Segre and Veronese embeddings.

Moreover, the previous formula  holds not only for embeddings~$f$,
but for any morphisms~$f$. As a consequence, one get
the desirable functoriality: \emph{if $f\colon Y\ra X$
is a morphism of projective algebraic varieties over~$F$
and $L\in\Pic(X)\otimes\R$, then $h_L\circ f\equiv h_{f^*L}$
(modulo bounded functions).}

Any function in the class~$h_L$ deserves to be called
a \emph{height function on~$X$ with respect to~$L$.}
However, it may be desirable to point out
specific height functions with good properties.
In the following paragraphs,
we show some cases where this is indeed possible.

\subsection{Algebraic dynamics (Tate, Silverman)}\label{ss.alg-dyn}
Let $X$ be a projective variety over~$F$, and assume
that it carries a dynamical system $\phi\colon X\ra X$,
and a real line bundle $L\in\Pic(X)\otimes\R$ such that
that $\phi^*L\simeq L^q$, for some real number $q>1$.
We shall say that $(X,\phi,L)$ is a \emph{polarized dynamical system},
and that $q$ is its weight.
Let $h_L^0$ be some arbitrary representative  of~$h_L$;
then, the following formula
\[ \hat h_L (x) = \lim_{n\ra\infty} q^{-n} h_L^0 (\phi^n(x)) \]
defines a height function $\hat h_L$ on~$X(\bar F)$
with respect to~$L$, which is independent of the choice of~$h_L^0$.
Moreover, it satisfies the following functional  equation
\[ \hat h_L (\phi(x)) = q \hat h_L(x), \quad\text{for any $x\in X(\bar F)$.}\]
In fact, it is the unique height function with respect
to~$L$ which satisfies this functional equation.
We cal it the \emph{canonical height function}.

Abelian varieties, that is, projective group varieties,
furnish especially beautiful examples of this situation.
If $X$ is an Abelian variety over~$F$, let $[n]$ be
the multiplication by an integer~$n$, an endomorphism of~$X$;
in particular, $[-1]$ is the inversion on~$X$.
Then, for any ample line bundle~$L$ on~$X$  which is symmetric
(that is, $[-1]^*L\simeq L$), one has
$[n]^*L\simeq L^{n^2}$ for any integer~$n$.
The various canonical height functions, for all integers~$n\geq2$,
coincide and are called the \emph{Néron-Tate height} on~$X$
relative to the line bundle~$L$.

Similarly, projective spaces, for the maps
$[x_0:\dots:x_n]\mapsto [x_0^q:\dots:x_n^q]$ (for some
integer~$q\geq 2$) and any line bundle,
and more generally toric varieties
are also interesting examples.

There are also nice examples for some K3-surfaces, first described
by~\cite{silverman1991}. (There, it is useful to work
with $\Pic(X)\otimes\R$, rather than $\Pic(X)$.)

\subsection{Height functions for geometric ground fields}
Assume that  $F=k(C)$ is the field of rational functions
on a regular curve~$C$ which is projective, geometrically
irreducible over a field~$k$.
Let $X$ be a projective variety over~$F$ and $L$
be a real line bundle on~$X$.
A projective $k$-variety $\mathscr X$ together
with a flat morphism $\pi\colon \mathscr X\ra C$,
the generic  fiber of which is~$X$, is called a \emph{model}
of~$X$ over~$C$; any line bundle $\mathscr L\in\Pic(\mathscr X)\otimes\R$
which gives back~$L$ on~$X$ is called a model of~$L$.

Now, let $x\in X(\bar F)$; it is defined
over a finite extension $E$ of~$F$ which is the field of rational
functions on a regular integral curve~$C'$, finite over~$C$.
By projectivity of~$X$ and regularity of~$C'$,
the point~$x$ is the generic fiber of a morphism $\eps_x\colon C'\ra X$.
Then, one can define
\[ h_{\mathscr L} (x) = \frac1{[C':C]} \deg_{C'} \eps_x^*\mathscr L. \]
The function $h_{\mathscr L}$ is a height function with
respect to~$L$.

\subsection{Arakelov geometry}
This point of view offers a sophisticated, and powerful,
way to mimic the geometric case
in order to obtain actual height functions in the arithmetic case.
Let $X$ be a projective variety over a number field~$F$,
let $L$ be a real line bundle on~$X$.

Let $\mathscr X$ be a model of~$X$ over the ring
of integers $\mathfrak o_F$, let $\mathscr L\in\Pic(\mathscr X)\otimes\R$
be a model of~$L$.
If we observe the analogy between function fields and number fields
under the point of view offered  by the product formula,
we see that $\mathfrak o_F$ behaves as the ring
of regular functions of an  affine curve. Consequently,
to get a height function,
we need to compactify somehow the spectrum of~$\mathfrak o_F$
taking into account the archimedean places of~$F$.
This is where Arakelov's ideas come in.

For any archimedean place~$v$ of~$F$, we set $\C_v=\overline{F_v}\simeq\C$
and indicate by a subscript~$v$ a $v$-adic completion,
or base change  to $v$-adic completion.
Let us endow the holomorphic line bundle~$L_v$ 
on the complex analytic space $X(\C_v)$ deduced from~$L$
with a continuous hermitian metric. Such a metric is a
way to define the \emph{size} of sections of~$L_v$.
It can be defined as a continuous function on the total space of~$L_v$ inducing
a hermitian norm on each fiber $L_{v,x}\simeq\C_v$
above a point $x\in X(\C_v)$. In other words,
any holomorphic section~$s$ of~$L_v$
over an open subset~$U$ of~$X(\C_v)$ is given a norm $\norm s$,
which is a continuous function $U\ra\R_+$, in such a way
that $\norm {fs}(x)=\abs{f(x)}\norm s(x)$ for any $x\in U$
and any holomorphic  function~$f$ on~$U$, and that
$\norm s(x)\neq 0$ if $s(x)\neq0$.

Let $\overline{\mathscr L}=(\mathscr L,(\norm{\cdot}_v))$
be the datum of such a model~$\mathscr L$ and
of hermitian metrics at all archimedean places~$v$ of~$F$.
It is called an hermitian line bundle over~$\mathscr X$.
(It is customary in Arakelov geometry to impose that
these metrics on~$\mathscr L$ be  conjugation invariant, but this hypothesis
is  not necessary here.)
Algebraic operations on line bundles such as taking duals,
or tensor products, can be done at the level of hermitian metrics,
so that there is a group $\hPic(\mathscr X)$
of isomorphism classes of hermitian line bundles on~$\mathscr X$.

Now, for $x\in X(\bar F)$, there is a finite extension~$E$
of~$F$ such that $x\in X(E)$,
and a morphism $\eps_x\colon\Spec \mathfrak o_{E}\ra\mathscr X$
which extends~$x$. We can then define
\[ h_{\overline{\mathscr L}}(x) =\frac1{[E:F]} \hdeg \eps_x^*\overline{\mathscr L}, \]
where $\hdeg$ means the Arakelov degree, an analogue for
hermitian line bundles over~$\Spec\mathfrak o_{E}$
of the geometric degree of line bundles  over complete curves.

Let us recall shortly the definition of this degree.
Let $\overline{\mathscr M}$ 
be a hermitian line bundle over~$\Spec\mathfrak o_E$. The
module of global sections is a projective $\mathfrak o_E$-module~$M$
of rank~1 and for any archimedean place~$v$ of~$E$, $M_v=M\otimes_E \C_v$ is endowed with a hermitian scalar product.
Then, one has
\[ \hdeg(\overline{\mathscr M}) = \log \frac{\Card( M/\mathfrak o_E m)}{\prod_v  \norm{m}_v }, \]
the right hand side being independent of the choice of a nonzero element $m\in M$. 
This independence follows 
from the fact that for any nonzero $a\in \mathfrak o_E$, the norm
of the ideal~$(a)$ coincides with the absolute value of the norm of~$a$.
In fact, it is an avatar of the product formula that was
used to define the height on projective spaces.

\subsection{Adelic metrics}\label{ss.adelic-metrics}
One may push the analogy between number fields and function fields
a bit further and do at non-archimedean places what Arakelov
geometry does at archimedean places.
This gives rise to the technique of adelic metrics, which
works both in the geometric and in the arithmetic settings.

Let $X$ be a projective variety over~$F$.
An adelic metric on~$L$ is a family $(\norm\cdot_v)$ of
continuous metrics on the line bundle~$L$  at all places~$v$
of~$F$ satisfying some ``adelic'' condition.

So let $v$ be a place of~$F$.
First, complete~$F$ for the absolute value given by~$v$,
then take its algebraic closure; this field admits
a unique absolute value extending~$v$; take its completion
for that absolute value.
Let $\C_v$ be the field ``of $v$-adic complex numbers''
so obtained; it is complete and algebraically closed.
A $v$-adic metric for the line bundle~$L$ 
can be defined similarly as in the case of archimedean
absolute values, as a continuous function on
the total space of the line bundle restricting to
the absolute value~$v$ on each fiber~$L_{x}\simeq \C_v$
at each point $x\in X(\C_v)$.

If $L$ is very ample and $(s_1,\dots,s_n)$ is a basis
of $H^0(X,L)$, then there is a unique metric on~$L$
such that
\[ \max(\norm{s_1}_v(x),\dots,\norm{s_n}_v(x))=1 \quad
\text{for all $x\in X(\C_v)$.} \]
Such a metric is called standard.
A family of  metrics on~$L$ will be called a standard
adelic metric if it is defined by this formula
for all places~$v$ of~$F$.

More generally,
a metric (or a family of metrics) on~$L$ will be called standard
if one can write $L\simeq L_1\otimes L_2^{-1}$ for two very ample
line bundles~$L_1$ and~$L_2$,
in such a way that the metric on~$L$ is the quotient
of standard metrics on~$L_1$ and~$L_2$.

However, the field~$\C_v$ is not locally compact,
so that  the resulting metrics lack good properties.
Therefore, one imposes the further condition that the
metric can be written as a standard metric times a
function of the form~$e^{\delta_v}$, where $\delta_v$
is a continuous and bounded function on~$X(\C_v)$.
Considering families of $v$-adic metrics, one imposes
that the function $\delta_v$ be identically~$0$ 
for all but finitely many places~$v$ of~$F$.

\section{Level sets for the canonical height}
We consider a polarized dynamical system $(X,L,\phi)$ over~$F$
with weight~$q>1$,
as in~\S\ref{ss.alg-dyn}.
Let $\hat h$ be its associated canonical height function,
satisfying the functional equation $\hat h(\phi(x))=q\hat h(x)$
for any $x\in X(\bar F)$.
The most important case will be the one associated
to Abelian varieties.

If $Y$ is a subvariety of~$X$ and $t$ is a real number,
we let $Y(t)$ be the set of points $x\in Y(\bar F)$
such that $\hat h(x)\leq t$.

\subsection{Preperiodic points}
Let $x\in X(\bar F)$. Its \emph{orbit} under~$\phi$
is the sequence of points of~$X$ obtained by iterating~$\phi$,
namely $(x,\phi(x),\phi^{(2)}(x),\dots)$.
One says that $x$ is periodic
if there exists $p\geq 1$ such that $x=\phi^{(p)}(x)$.
One says that $x$ is preperiodic if its orbit is finite
or, equivalently, if there are integers~$n\geq 0$ and~$p\geq 1$
such that $\phi^{(n)}(x)=\phi^{(n+p)}(x)$.

When $X$ is an Abelian variety and $\phi$ is the multiplication
by an integer~$d\geq 2$, preperiodic points are exactly
the torsion points of~$X$. One direction is clear: if $x$ has
finite order, say~$m$, then every multiple of~$x$ is killed by
the multiplication by~$m$; since there are finitely many points~$a\in X(\bar F)$
such that $[m]a=0$, the orbit of~$x$ is finite.
Conversely, if the orbit of~$x$ is finite, let $n$ and~$p\geq 1$
be integers such that $\phi^{(n)}(x)=\phi^{(n+p)}(x)$;
this implies $[d^n]x=[d^{n+p}]x$, hence $[d^n(d^p-1)]x=0$,
so that $x$ is a torsion point.

The canonical height of a preperiodic point must be zero. Indeed,
let $x$ be preperiodic and let $n$ and~$p\geq 1$ be integers
such that $\phi^{(n)}(x)=\phi^{(n+p)}(x)$. Computing
the canonical height of both sides of the equality,
we get $q^n \hat h(x)=q^{n+p}\hat h(x)$, hence $\hat h(x)=0$
since $q>1$, so that $q^n\neq q^{n+p}$.

\subsection{Points of canonical height zero}
If $F$ is a number field, or a function field over a finite field,
the converse holds. Let $x\in X(\bar F)$ be such that $\hat h(x)=0$.
Let $E$ be a finite extension of~$F$ such that $x\in X(E)$.
Any point in the orbit of~$x$ has canonical height zero;
by Northcott's finiteness theorem, the orbit of~$x$ is finite.
In fact, this statement was the main result of~\cite{northcott1950}!

However, if $F$ is  a function field over an algebraically closed
field~$k$, Northcott's finiteness theorem is false
and this property does not hold anymore.
Indeed, we can consider a ``constant'' dynamical system $(Y,M,\psi)$
defined over~$k$ and view it over~$F$.
Then, all points in $Y(\bar k)$ have canonical height zero,
unless $k$ is the algebraic
closure of a finite field,
they are usually not preperiodic. 

In fact, a theorem of~\cite{chatzidakis-hrushovski2008}
shows that this obstruction is essentially the only one.
This generalizes an old result of~\cite{lang-neron1959}
for Abelian varieties.
Because it is  simpler to quote, let us only
give the particular case due to~\cite{baker2009}.
\emph{Let  $X=\P^1$ and $\phi\in F(t)$ be a rational function
of degree~$q\geq 2$, then there exists a non-preperiodic
point $x\in X(\bar F)$ such that $\hat h(x)=0$
if and only $\phi$ is conjugate (by a homography in~$\bar F$)
to a rational function~$\psi\in k(t)$.}

\subsection{The geometry of points of canonical height zero}
In the 60s, motivated by the conjecture
of Mordell and its extension by Lang,
Manin and Mumford  had conjectured 
that an integral subvariety of a complex Abelian variety
cannot contain a dense set of torsion points, 
unless the subvariety is the translate of an Abelian
subvariety by a torsion point.

This expectation was proved to be a theorem, due to~\cite{raynaud83c}
(there are many other proofs now). 
We quote it in the slightly different, but equivalent, form:
\emph{Assume that $X$
is an Abelian variety over an algebraically closed field
of characteristic zero and let $Y$
be a closed subvariety of~$X$. Then the Zariski closure of~$Y(0)$
is a finite union of translates of Abelian  subvarieties by
torsion points.}
In fact, the proof relies on techniques from arithmetic geometry,
and its crucial part assumes that 
$X$ and $Y$ are defined over a number field.

Similarly, the study of the analogue of Manin \&\ Mumford's question 
over algebraically closed fields of positive characteristic 
would reduce to the case of function fields.
But there, eventual constant Abelian varieties within~$X$ 
create difficulties. The precise theorem has been 
first proved by~\cite{scanlon2001,scanlon2005}. Before I state it, 
let me recall the existence of the \emph{Chow trace} of~$X$,
a ``maximal''  Abelian variety~$X'$ defined over~$k$ together
with a morphism $X'\otimes_k F\ra X$.
Then, \emph{if $Y$ is a subvariety of~$X$,
the Zariski closure of~$Y(0)$
is a finite union of varieties~$Z$ such that the quotient
of~$Z$ by its stabilizer~$G_Z$ is a translate
of a subvariety of~$(X/G_Z)'\otimes_k F$ defined over~$k$
by a torsion point of~$X/G_Z$.}

\subsection{}
In the context of dynamical systems, the
question of Manin \&\ Mumford generalizes as follows.
Let $Y$ be a subvariety of~$X$ and let $Y(0)$
be the set of points $x\in X(\bar K)$ such that $\hat h(x)=0$.
Is it true that $Y(0)$ is not dense in~$Y$, unless  this
is somewhat explained by the geometry of~$Y$ with respect to~$\phi$,
for example, unless $Y$ is itself preperiodic?
However, the answer to this basic expectation was proved to be \emph{false},
by an example of Ghioca \&\ Tucker. A subsequent
paper by~\citet*{ghioca-tucker-zhang2011} tries to
correct the basic prediction.

\subsection{The conjecture of Bogomolov}
Still in conjunction with Mordell's conjecture,
\cite{bogomolov80b} had strengthtend 
Manin--Mumford's question by requiring to prove that over number
field, if $C$ is a curve over genus~$\geq 2$  embedded in its Jacobian~$J$,
there exists a positive real number~$\eps$
such that $C(\eps)=\{x\in C(\bar F)\,;\, \hat h(x)\leq\eps\}$
is finite.

This conjecture has been generalized by~\cite{zhang95b}
to subvarieties of Abelian varieties over a number field: if $Y$
is a subvariety of an Abelian variety~$X$, does there exists
a positive real number~$\eps$ such that $Y(\eps)$
is contained in a finite union of translates of Abelian
subvarieties of~$X$ by torsion points, contained in~$Y$?  In other
words, is it true that $Y(\eps)\subset \overline{Y(0)}$
for small enough $\eps$?

These two questions have been solved positively
by~\cite{ullmo98} and~\cite{zhang98} respectively.
They make a crucial use of equidistribution arguments
that will be explained below.
Soon after, \cite{david-p98} gave another proof.

The analogous case of function fields is mostly open, the last part
of this text will be devoted to explaining
how Gubler had been able to use
ideas of equidistribution to prove important cases in this setting.
Note that over a function field of characteristic zero,
the case of a curve embedded in its Jacobian 
(Bogomolov's original on) has been settled positively 
by \cite{faber2009b} and \cite{cinkir2011} using 
a formula of~\cite{zhang2010} for the height of the Gross-Schoen cycle.

\section{Equidistribution (arithmetic case)}

\subsection{}
Equidistribution is a prevalent theme of analytic number theory:
it is a (partially) quantitative way of describing how
discrete objects collectively model a continuous phenomenon.
The most famous result is probably the equidistribution
modulo~$1$ of multiples $n\alpha$ of some fixed irrational number~$\alpha$,
due to Weyl.

Here, we are interested in algebraic points~$x\in X(\bar F)$ of a variety~$X$
defined over~$F$.
To have some chance of getting some continuous phenomenon,
we consider, not only the points themselves,
but also their conjugates, that is, the full orbit of those points under
the Galois group~$\Gal(\bar F/F)$.
The continuous phenomenon requires some topology,
so we fix a place~$v$ of~$F$ and an embedding of~$\bar F$
into the field~$\C_v$.

Let $x$ be any point in $X(\bar F)$.
Viewed from the field~$F$, the point~$x$ is not discernible
from any of its conjugates  $x_1=x,\dots,x_m$ which are
obtained from~$x$ by letting the group of $F$-automorphisms of~$\bar F$
act. So we define a probability measure $\mu(x)$ on~$X(\C_v)$ by
\[ \mu(x) = \frac1m \sum_{j=1}^m \delta_{x_j}, \]
where $\delta_{x_j}$ is the Dirac measure at the point~$x_j\in X(\bar F)\subset X(\C_v)$.

The first equidistribution result in this field
is the following.

\begin{theo}[\citet*{szpiro-u-z97}]\label{theo.suz}
Assume that $F$ is a number field and $v$ is an archimedean place of~$F$.
Let $X$ be an Abelian variety over~$F$ and
let $(x_n)$ be a sequence of points in~$X(\bar F)$
satisfying the following two assumptions:
\begin{itemize}
\item The Néron--Tate height of~$x_n$ goes to~$0$ when $n\ra\infty$;
\item For any subvariety~$Y$ of~$X$ such that $Y\neq X$, the
set of indices~$n$ such that $x_n\in Y$ is finite.
\end{itemize}
Then the sequence of probability measures $(\mu(x_n))$
on the complex torus $X(\C_v)$ converges vaguely to the normalized
Haar measure of $X(\C_v)$.
\end{theo}

The proof uses Arakelov geometry and holds in a wider context
than that of Abelian varieties. We shall see more about it
shortly but I would like to describe how \cite{ullmo98}
and \cite{zhang98} used those ideas to obtain a \emph{proof}
of Bogomolov's conjecture.

\subsection{}
So assume that $Y$ is a subvariety of~$X$, with $Y\neq X$, containing
a sequence $(x_n)$ of algebraic points such that $\hat h(x_n)\ra 0$
and which is dense in~$Y$ for the Zariski topology.
We want to show that $Y$ is a translate of an Abelian subvariety
of~$X$ by a torsion point. To that aim, we may mod out~$X$ and~$Y$
by the stabilizer of~$Y$.
The definition of the Néron--Tate height on~$X$ comes
from some ample line bundle; its Riemann form on~$X(\C_v)$
is a positive differential form~$\omega$ of bidegree~$(1,1)$.

Now, a geometric result implies that there exists a positive
integer~$m$ such that the map
\[ \phi \colon Y^m\ra X^{m-1}, \quad (y_1,\dots,y_m)\mapsto
(y_2-y_1,\dots,y_m-y_{m-1}) \]
is generically finite.
From the sequence~$(x_n)$, one constructs a similar sequence $(y_n)$
of points in $Y^m(\bar F)$ whose heights converge to zero
and which are Zariski dense in~$Y^m$;
more precisely, for any  subvariety~$Z$ of~$Y^m$
such that $Z\neq Y^m$,
the set of indices~$n$ such that $y_n\in Z(\bar F)$ is finite.
A variant of Theorem~\ref{theo.suz} (see also Theorem~\ref{theo.yuan} below)
implies that
the sequence $\mu(y_n)$ of probability measures
converges to the canonical probability measure on~$Y^m(\C_v)$
given by the differential form
$(\omega_1+\dots+\omega_m)^{md}$ on the smooth locus of~$Y^m$.
(Here, $d=\dim (Y)$ and $\omega_j$ means the differential form on~$Y^m$
coming from~$\omega$ on the $j$th factor of~$Y^m$.)
Write $\mu(Y^m)$ for this measure; in fact, one has $\mu(Y^m)=\mu(Y)^m$.
So we have the equidistribution property
\[ \mu(y_n) \ra \mu(Y)^m. \]
If we apply the map~$\phi$, we get automatically
\[ \mu(\phi(y_n))\ra \phi_* \mu(Y)^m. \]

On the other hand, the sequence $(\phi(y_n))$ also satisfies
an equidistribution property, but the limit measure being
$\mu(\phi(Y^m))$. This implies  an \emph{equality} of probability measures
\[ \phi_*\mu(Y)^m=\mu(\phi(Y^m))) ,\]
a geometric refinement of the initial fact
that $\phi$ is generically finite with image $\phi(Y^m)$.

However, both sides of this equality come from differential forms,
and this equality implies that the differential forms
$ (\omega_1+\dots+\omega_m)^{md} $
and $\phi^* (\omega_1+\dots+\omega_{m-1})^{md}$ on~$Y^m$
coincide up to a constant multiple.

The contradiction comes from the fact that
$(\omega_1+\dots+\omega_m)$ is strictly positive everywhere
(at least, on the smooth locus of~$Y^m$)
while $\phi^*(\omega_1+\dots+\omega_{m-1})^{md}$ vanishes
where $\phi$ is not étale (that is, not a local diffeomorphism),
in particular on the diagonal of~$Y^m$.
To be resolved, this contradiction requires that $md=0$, hence
that $Y$ is a point, necessarily a torsion point.

\subsection{Heights for subvarieties}
One of the major ingredients in the proof of the
equidistribution theorem is a notion of a height
with respect to a metrized line bundle,
not only for points, but for all subvarieties.
The definition, first introduced by~\cite{faltings91},
goes as follows.

We assume that $X$ is a projective variety over a number field~$F$
and $L$ is a line bundle on~$X$.
Let $\mathscr X$ be a projective flat model over
the ring~$\mathfrak o_F$ and $\mathscr L$
be line bundle on~$\mathscr X$ which is a model of~$L$.
We also assume that $L$ is endowed 
with smooth hermitian metrics at all archimedean places of~$F$.

From these metrics, complex differential  geometry defines
differential forms  $c_1(\overline L_v)$ on the complex analytic
varieties $X(\C_v)$, for all archimedean places~$v$ of~$F$.
This form is called the \emph{first Chern form},
or the \emph{curvature form},
of the hermitian line bundle~$\bar L_v$;
it is a representative of the first Chern class of~$L$
in the De Rham cohomology of~$X(\C_v)$.
It is really a fundamental tool in complex algebraic geometry.
For example, when $X$ is smooth, say,
the Kodaira embedding theorem asserts that $L$ is ample
if and only if it possesses a hermitian metric such that
its curvature form is positive definite on each tangent
space of~$X$.
I refer to Section~1.4 of~\cite{griffiths-h78} for more details.

Faltings's definition of the height of an irreducible closed
subvariety $\mathscr Y\subset\mathscr X$ is by induction
on its dimension.

If $\dim(\mathscr Y)=0$, then $\mathscr Y$ is a closed point;
then, its residue field $\kappa(\mathscr Y)$ is a finite field
and one defines
\begin{equation}
h_{\overline{\mathscr L}}(\mathscr Y)= \log \Card(\kappa(\mathscr Y)).
\end{equation}
Otherwise, one can consider a (nonzero) meromorphic section~$s$ of
some power~$\mathscr L^m$ of~$\mathscr L$  on~$\mathscr Y$.
Its divisor $\div(s)$
is a formal linear combination  of irreducible
closed subschemes $\mathscr Z_j$
of~$\mathscr Y$, with multiplicities~$a_j$ (the order of vanishing,
or minus the order of the pole of~$s$ along~$\mathscr Z_j$)
and
\begin{equation}
h_{\overline{\mathscr L}}(\mathscr Y)
= \frac1m \sum a_j h_{\overline{\mathscr L}}(\mathscr Z_j)
+ \sum_{v} \int_{Y(\C_v)} \log\norm s^{-1/m} c_1(\overline L_v)^{\dim Y}
\end{equation}
where $Y=\mathscr Y\otimes F$ and $v$ runs over archimedean places of~$F$.
In fact, the right hand side of this formula does not depend
on the choice of~$s$.

One can prove that this new definition recovers the previous one for points.
More precisely, let $y\in X(\bar F)$, let $Y\in X$ be the corresponding
closed point  and let $\mathscr Y$ be its Zariski closure in~$\mathscr X$.
Then,
\[ h_{\overline{\mathscr L}}(\mathscr Y) = \deg(Y) h(y), \]
where $\deg(Y)$ is the degree of the closed point~$Y$,
or the degree of~$Y$ as a subvariety of~$X$ with respect to the line bundle~$L$.

\begin{prop}[\cite{zhang95b}]
\label{prop.zhang}
Assume that $L$ is ample, that $\mathscr L$ is relatively
numerically effective and that the curvature forms
$c_1(\overline L_v)$ are nonnegative
for any archimedean place~$v$ of~$F$.

Let $(x_n)$ be a sequence of points in $X(\bar F)$.
Assume that for any subvariety~$Y$ of~$X$ such that $Y\neq X$, the set
of indices~$n$ such that $x_n\in Y$ is finite.
Then,
\begin{equation}\label{eq.zhang}
 \liminf_n h_{\overline{\mathscr L}}(x_n) \geq \frac{ h_{\overline{\mathscr L}}(\mathscr X)}{\dim(\mathscr X) \deg_L(X)}.\end{equation}
\end{prop}

This proposition follows easily from a (difficult) theorem in Arakelov
geometry that implies the existence of global sections over~$\mathscr X$
of large powers $\mathscr L^{m}$ which have controlled norms.
Using those sections in the inductive definition of the height
leads readily to the indicated inequality.

In presence of a sequence~$(x_n)$ such 
that $h_{\overline{\mathscr L}}(x_n)$ converges
to the right hand side of Inequality~\ref{eq.zhang},
\citet*{szpiro-u-z97} proved that the probability measures
$\mu(x_n)$ equidistribute towards the measure
$\mu_{X,v}= c_1(\bar L_v)^{\dim(X)}/\deg_L(X)$ on $X(\C_v)$.
The heart of the proof is to apply the fundamental inequality~\eqref{eq.zhang}
for small perturbations of the hermitian metrics,
as a \emph{variational principle.}
Since $X(\C_v)$ is compact and metrizable, the space  of probability measures
on $X(\C_v)$ is metrizable and compact, so we may assume
that $\mu(x_n)$ converges to some limit~$\mu$ and we need
to prove that $\mu$ is proportional to $c_1(\bar L_v)^{\dim X}$.

Let us multiply the metric on~$\bar L_v$ by some function
of the form $e^{-\eps\phi}$, where $\phi$ is a smooth function
on $X(\C_v)$. Then, the left hand side of the inequality~\eqref{eq.zhang}
becomes
\[ \lim_n \left( h_{\overline{\mathscr L}}(x_n) + \eps \int_{X(\C_v)}\phi\,\mathrm d\mu(x_n)\right)
=  \frac{ h_{\overline{\mathscr L}}(\mathscr X)}{\dim(\mathscr X) \deg_L(X)}
+ \eps \int_{X(\C_v)} \phi\,\mathrm d\mu, \]
while its right hand side
is
\[  \frac{ h_{\overline{\mathscr L}}(\mathscr X)}{\dim(\mathscr X) \deg_L(X)}
+ \eps \int_{X(\C_v)} \phi \frac{c_1(\bar L_v)^{\dim X}}{\deg_L(X)}
+ \mathrm O(\eps^2). \]
Consequently, when $\eps\ra 0$,
\[ \eps (  \mu(\phi)  - \mu_{X,v}(\phi) ) \geq \mathrm O(\eps^2). \]
For small positive $\eps$, we get $\mu(\phi)\geq\mu_{X,v}(\phi)$,
and we have the opposite inequality for small negative $\eps$.
Consequently, $\mu(\phi)=\mu_{X,v}(\phi)$, hence the equality $\mu=\mu_{X,v}$.

A subtle point of the proof lies in the possibility
of applying Proposition~\ref{prop.zhang}
to the modified line bundle.
When the curvature form $c_1(\overline{L})$ is strictly positive,
then it remains so for small perturbations,
hence the proof is legitimate. This is what happens
in Theorem~\ref{theo.suz}, and what is needed
for the proof  of Bogomolov's conjecture by Ullmo and Zhang.

Inspired by an inequality of Siu and the holomorphic
Morse inequalities of Demailly, \cite{yuan2008} proved
the following  general equidistribution theorem.

\begin{theo}[\cite{yuan2008}]
\label{theo.yuan}
Assume that $F$ is a number field and $v$ is an archimedean place of~$F$.
Let $X$ be an algebraic  variety over~$F$, let $\bar L$
be an ample line bundle on~$X$ with a semi-positive adelic metric.
Let $(x_n)$ be a sequence of points in~$X(\bar F)$
satisfying the following two assumptions:
\begin{itemize}
\item The heights of~$x_n$ with respect to~$\bar L$
converge to~$0$ when $n\ra\infty$;
\item For any subvariety~$Y$ of~$X$ such that $Y\neq X$, the
set of indices~$n$ such that $x_n\in Y$ is finite.
\end{itemize}
Then the sequence of probability measures $(\mu(x_n))$
on the complex space $X(\C_v)$ converges to the
unique probability measure proportional to $c_1(\bar L)^{\dim(X)}$.
\end{theo}
I cannot say much more on this here, and I
must refer the reader to the paper of~\cite{yuan2008}.

Observe anyway that under the indicated hypotheses, $c_1(\bar L)$
is not necessarily a differential form, but only a positive
current of bidegree~$(1,1)$. Consequently, the definition of the measure
$c_1(\bar L)^{\dim X}$ requires some work.
It goes back to fundamental work in pluripotential
theory by \cite{bedford-t82} and~\cite{demailly1985}.
In our setting, it can be defined by an approximation process,
considering sequences of smooth \emph{positive} hermitian metrics on~$L$
which converge uniformly to the initial metric.
See my survey~\cite{chambert-loir2011} for more details.

\section{Measures on analytic spaces}

\subsection{}
Our setting is that of a global field~$F$. Let $X$ be
a projective algebraic (irreducible) variety over~$F$.
For any place~$v$ of~$F$, we will consider the analytic
space~$X_v^\an$ associated to~$v$.

If $v$ is archimedean, then $X_v^\an=X(\C_v)$
is the set of complex points of~$X$, where $\C_v=\C$
is viewed as an $F$-algebra via the embedding corresponding
to the place~$v$.

When $v$ is non-archimedean, then $X_v^\an$
is the analytic space over the complete algebraically closed
field~$\C_v$,
as defined by~\cite{berkovich1990}.
I must refer to the other contributions in this volume for
background on Berkovich spaces, as well
as to those of~\cite{baker2008b} and~\cite{conrad2008}
in the proceedings of the 2007~Arizona conference
edited by~\citeauthor{zz-arizona2007}.
Here, I will content myself
with the few  following comments. First of all, $X_v^\an$
is a reasonable topological space: it is compact
and  locally pathwise connected.
It is even locally contractible;
this is the main theorem of \cite{berkovich1999} when $X$ is smooth,
recently extended to the general case by \cite{hrushovski-loeser2009}.
Moreover, $X_v^\an$ contains the set $X(\C_v)$
as a dense subset, and the topology of $X_v^\an$ restricts to its
natural (totally disconnected) topology on $X(\C_v)$.
So $X_v^\an$ has many other points than those of~$X(\C_v)$,
some which will play a crucial role below.

\subsection{}
To fix ideas, assume that we are in the geometric case,
so that $F=k(C)$ is the field of functions on  a curve~$C$.
Let $\mathscr X$ be a projective model of~$X$ over the curve~$C$
and let $\mathscr L$ be a line bundle on~$\mathscr X$;
let $L$ be its restriction to~$X$.

For any place~$v$ of~$F$, the model~$\mathscr L$
gives rise to a {\og $v$-adic metric on~$L$\fg}.
This notion is similar to what was discuted in the case of~$X(\C_v)$; 
in particular, any section~$s$ of~$L$ on an open subset~$U$ of~$X$
has a norm $\norm s$ which is a continuous function on the
corresponding subset~$U_v^\an$ of~$X_v^\an$, and
does not vanish if $s$ does not vanish.
The construction  of the metric is also similar to that
of standard metrics.
Assume first that  $\mathscr L$ be very ample; then,
the metric on~$L$ is the unique  metric such that 
for any  generating set $(s_1,\dots,s_n)$
of the module $\Gamma(\mathscr X,\mathscr L)$
of integral global sections, one has
\[ \max(\norm{s_1}(x),\dots,\norm{s_n}(x))=1 \quad
\text{for any $x\in X_v^\an$.} \]
In general, one can at least write
$\mathscr L$ as the difference~$\mathscr L_1\otimes\mathscr L_2^{-1}$
of two very ample line bundles on~$\mathscr X$,
and the metric on~$L$ is the quotient of the to metrics
given by the models~$\mathscr L_1$ and~$\mathscr L_2$.

I claim that there exists a measure,
written $c_1(\overline L)_v^{\dim X}$, on $X_v^\an$
such that for any nonzero global section of~$L$,
\[ h_{\mathscr L} (\mathscr X) = h_{\mathscr L}(\overline{\div(s)})
+ \sum_v \int_{X_v^\an} \log\norm{s}_v^{-1} c_1(\bar L)_v^{\dim X}, \]
where $\overline{\div(s)}$ is the Zariski closure in~$\mathscr X$
of the divisor of~$s$.
More generally, if $Y$ is an integral subvariety of~$X$,
with Zariski closure~$\mathscr Y$ in~$\mathscr X$,
and if $s$ is any nonzero global section of~$L|_Y$, 
one can define a measure $c_1(\overline L)_v^{\dim Y}\delta_Y$
on $Y_v^\an$ such that
\[ h_{\mathscr L} (\mathscr Y) = h_{\mathscr L}(\overline{\div(s)})
+ \sum_v \int_{X_v^\an} \log\norm{s}_v^{-1} c_1(\bar L)_v^{\dim Y}\delta_Y. \]

This measure is defined as follows,
see \cite{chambert-loir2006}; the presentation given
here, using algebraically closed valued fields,
is due to \cite{gubler2007a}.

The analytic space~$\mathscr X_v$
admits a canonical reduction $\mathsf X_v$ over
the residue field of~$\C_v$,
which maps to the natural reduction of~$\mathscr X$.
Moreover, there is a reduction map $X_v^\an\ra \mathsf X_v$
and
the generic point of an irreducible component~$\mathsf Z$ of~$\mathsf X_v$ is
the image of a unique point~$z$ of~$X_v^\an$.
By functoriality, one also has a line bundle $\mathsf L_v$ on~$\mathscr X_v$.
The measure $c_1(\bar L)^{\dim X}_v$ is the following
linear combination of Dirac measures:
\[ c_1(\bar L)^{\dim X}_v = \sum_{\mathsf Z}   (c_1(\mathsf L_v)^{\dim X}|\mathsf Z) \delta_z \]
where the coefficients $(c_1(\mathsf L_v)^{\dim X}|\mathsf Z)$
are given by usual (numerical) intersection theory
and the sum is over the irreducible components~$\mathsf Z$
of~$\mathsf X_v$.

This measure is positive if $\mathscr L$ is relatively
ample, and its total mass is equal to the degree of~$X$ with respect
to~$L$.

The definition of the measure $c_1(\bar L)_v^{\dim Y}\delta_Y$ is analogous.

Up to its measure-theoretic formulation,
the validity of the asserted formula for heights follows
from work of~\cite{gubler1998}.

\subsection{}
\cite{zhang95b} had defined a notion of semipositive metrics,
which are defined as uniform limits of metrics given by
models $(\mathscr X,\mathscr L)$, where $\mathscr L$
is relatively numerically effective---that is, gives a nonnegative
degree to any vertical subvariety.
He also showed that semipositive metrized line bundles
allow to define heights of subvarieties by approximation
from the case of models/classical Arakelov geometry.

Adapting this construction I defined in~\cite{chambert-loir2006}
the measures $c_1(\bar L)^{\dim X}_v$ by approximation
from the above definition in the case of models.
In the end, the proof is very close to that of the existence
of products of positive $(1,1)$-currents by~\cite{bedford-t82}.
(In fact, this article
only considers projective varieties over a local $p$-adic field;
the general case has been treated by~\cite{gubler2007a},
in a similar fashion.)

\subsection{}
As we have shown in~\cite{chambert-loir-thuillier2009},
these measures can be used to recover the heights defined
by~\cite{zhang95b}.
Namely, if  $\bar L$ is a line bundle
on~$X$ with a semi-positive adelic metric,
$Y$ is an integral subvariety of~$X$, Zhang defined 
the height $h_{\bar L}(Y)$ of~$Y$ with respect to~$\bar L$.
For any regular meromorphic section~$s$ of $L|_Y$, one has
\[ h_{\bar L} (Y) = h_{\bar L}({\div(s)})
+ \sum_v \int_{X_v^\an} \log\norm{s}_v^{-1} c_1(\bar L)_v^{\dim Y}\delta_Y. \]

In the case of curves, and a few
cases in higher dimensions, I showed in~\cite{chambert-loir2006}
that they also give rise to equidistribution theorems
totally analogous  to the one of~\citet*{szpiro-u-z97}.
The article of~\cite{yuan2008} proved what can be considered
the most general equidistribution theorem possible in this context.
Namely, the non-archimedean analogue of Theorem~\ref{theo.yuan}
still holds. While that paper restricts to the case
of number fields, its ideas have been transposed to the case 
of function fields by \cite{faber2009} and~\cite{gubler2008}.


I also refer to~\cite{petsche2010} for a general discussion
of convergence of measures on the Berkovich projective space,
as well as for a non-archimedean
analogue of Weyl's equidistribution criterion.

\subsection{}
In some cases, one can deduce from these equidistribution theorems
explicit results in algebraic number theory.
Let us give an example in the case of the line bundle~$L=\mathscr O(1)$
on  the projective line $X=\P^1$, with its metrization giving rise
to the standard height. Fix an ultrametric place~$v$ of~$F$.
Then, the measure  $c_1(\bar L)_v$  on~$X_v^\an$ is the Dirac measure
at a particular point~$\gamma$, called the Gauss-point
because it corresponds to the Gauss-norm on the
algebra $F_v[T]$ (viewed as the algebra of functions on
the affine line~$\A^1=\P^1\setminus\{\infty\}$).
So in this case, the equidistribution theorem asserts that
for any sequence $(x_n)$ of distinct points on~$X(\bar F)$
such that $h(x_n)\ra 0$, the measures $\mu(x_n)$ on~$X^\an_v$ converge
to the Dirac measure $\delta_{\gamma}$.

This gives a strong constraint on such sequences.
For example, it is impossible that all~$x_n$ 
be totally $v$-adic (an algebraic point is ``totally $v$-adic''
if all of its conjugates are defined over~$F_v$). 
Indeed, if $x_n$ is totally $v$-adic,
then the measure $\mu(x_n)$ is supported by the compact
subset~$X(F_v)$ of~$X^\an_v$.
If all $x_n$ were totally $v$-adic, the limit measure of~$\mu(x_n)$
would be supported by~$X(F_v)$, but the Gauss-point  does \emph{not}
belong to~$X(F_v)$. Similar results were proved by~\cite{baker-hsia2005}.

\section{Bogomolov's conjecture for totally degenerate abelian varieties}

\subsection{}
\cite{gubler2007b} had the idea of using these measures
to attack the unsolved Bogomolov  conjecture over function fields,
using equidistribution theorems for points of small height
at some place of the ground field to get a proof of
the conjecture following the strategy of~\cite{ullmo98,zhang98}.

So let $F$ be a function field and let $v$ be a place of~$F$.
Let $X$ be an Abelian variety over~$F$, let $\bar L$ be an ample symmetric
line bundle on~$X$ with its canonical adelic metric that
gives rise to the Néron-Tate height~$\hat h$.
Let $Y$ be a closed integral subvariety of~$X$
which is not the translate of 
an Abelian subvariety by a torsion point.
We want to prove that for some positive~$\eps$,
$Y(\eps)$ is not Zariski-dense. Assume the contrary.
We would then want to construct a sequence~$(y_n)$
in~$Y(\bar F)$ satisfying the assumptions of the equidistribution
theorem, namely $\hat h(y_n)\ra 0$ and for any subvariety~$Z$
of~$Y$ such that $Z\neq Y$, the set of integers~$n$ such that $y_n\in Z$
is finite. However, the set of subvarieties of~$Y$ may be uncountable,
hence such a sequence may not exist.
Anyway, one can construct a \emph{net} $(y_n)$, that is 
a family of points indexed by a filtered ordered set~$N$,
such that $h(y_n)\ra 0$ and for any subvariety~$Z\subsetneq Y$,
the set of indices~$n$ such that $y_n\in Z$ is bounded in~$N$.
The statement and the proof of the equidistribution theorem
adapt readily to this case.

We redo the same geometric reduction, assuming that
the stabilizer of~$Y$ is trivial, and that the morphism
$\phi\colon Y^m\ra X^{m-1}$ given by $(y_1,\dots,u_m)\mapsto (y_2-y_1,\dots,y_m-y_{m-1})$ is generically finite, with image~$Z$.
As above, we construct a generic net~$(y_n)$
of small points in~$Y^m$ whose image $(\phi(y_n))$ s a
generic net of small points in~$Z$. 
This gives two equidistribution
theorems in the Berkovich spaces~$(Y^m)^\an_v$ and~$Z^\an_v$ 
at the chosen place~$v$,
with respect to canonical measures 
$\mu_v(Y^m)=c_1(\bar L|_{Y^m})_v^{m\dim Y}$
and $\mu_v(Z)=c_1(\bar L|_Z)^{\dim Z}$,
where we write $\bar L|_{Y^m}$ and $\bar L|_Z$ 
for the metrized line bundles
on~$Y^m$ and~$Z$ deduced from those naturally given by~$\bar L$
on~$X^m$ and~$X^{m-1}$.
By construction, $\phi_*\mu_v(Y^m)=\mu_v(Z)$.

To get a contradiction, we need to have more information about
these measures.

\subsection{}
If $X$ has good reduction at~$v$, the very definition
of the measure $\mu_v(X)$ shows that it is the Dirac measure
at a single point of~$X^\an_v$. Indeed, let $\mathscr X$ be the Néron model of~$X$
over the ring of integers~$\mathfrak o_v$ of~$F_v$; since $X$
has good reduction, $\mathscr X$ is proper and smooth,
and its special fiber is an Abelian variety. Then, one can show
that the generic point of this fiber has a unique preimage~$\xi$
under the reduction map from the Berkovich space~$X^\an_v$ to the special fiber.
One has $\mu_v(X)=\deg_L(X) \delta_{\xi}$.

In the case where all of~$X$, $Y$ and~$Z$ have good reduction at~$v$
(this happens for almost all places~$v$), the measures
$\mu_v(Y^m)$ and $\mu_v(Z)$ are supported at a single point
and the equality of measures $\phi_*\mu_v(Y^m)=\mu_v(Z)$
gives no contradiction.

Also, if $X$ has good reduction, the measures~$\mu_v(Y^m)$
and~$\mu_v(Z)$ will be supported at finitely many points
and it will still be difficult to draw a contradiction.

\subsection{}
Consequently, to succeed, this equidistribution approach
needs to consider places of bad reduction of~$X$.
The case treated by~\cite{gubler2007b}
is the one of  \emph{totally degenerate Abelian varieties},
those being as far as possible from Abelian varieties of good reduction.
Recall that any Abelian variety over~$F$
has a canonical model over the local ring $\mathfrak o_{F,v}$
at the place~$v$, called its Néron model. Possibly
after a finite extension of the ground field,
the connected
component of the identity of the special fiber
of the Néron model is an extension of torus~$\gm^a$
by an Abelian variety; totally degenerate
Abelian varieties are those for which this torus has dimension~$\dim(X)$.

Assume this is the case and set $g=\dim(X)$.
Possibly after some finite extension of~$F$,
By theorems of Tate, Raynaud, Bosch, Lütkebohmert
in Tate's setting of rigid analytic spaces,
extended to the Berkovich context in~\cite[\S6.5]{berkovich1990},
the analytic space~$X^\an_v$ associated
to the Abelian variety~$X$ can be written as the quotient of a torus $\gmvg$
by a discrete subgroup~$\Omega$ of rank~$g$ in~$\gm^g(F_v)$.
In fact,  the torus~$\gmvg$ is the universal cover
of the Berkovich space~$X^\an_v$.

In particular, the topological fundamental group
of the analytic space associated
to our totally degenerate Abelian variety~$X$ is isomorphic to~$\Z^g$.
This does not reflect however the richness of étale covers
of Abelian varieties ---  the fundamental group of a complex Abelian
varieties of dimension~$g$ his~$\Z^{2g}$, while 
the $\ell$-adic fundamental group would be $\Z_\ell^{2g}$ (provided
$\ell$ is distinct from the characteristic of the ground field).
This indicates that, in some sense, 
the reduction at archimedean places is at least twice as bad 
as the worst possible ultrametric places of bad reduction.

Here enters tropical geometry.
\subsection{}
We first analyze the tropicalization of a torus.
By definition, the Berkovich space of~$\gm$ at the place~$v$
is the set of all multiplicative seminorms on the ring~$F_v[T,T^{-1}]$
which extend the fixed absolute value on~$F_v$.
So there is a natural map from~$\gmv$ to the real line~$\R$
that maps a semi-norm~$\chi$ to the real number $-\log\abs{\chi(T)}$.
In fact, the semi-norm~$\chi$ is viewed as a point~$x$ of~$\gmv$,
and $\abs{\chi(T)}$ is viewed as $\abs{T(x)}$, so that a more natural
way to write this map is $\tau\colon x\mapsto -\log\abs{T(x)}$.
An even more natural way would be to consider the
map $x\mapsto \abs{T(x)}$ from~$\gmv$ to~$\R_+^*$, because it
does not require the choice of a logarithm function.

This ``tropicalization'' map~$\tau$ is continuous and surjective.
It has a canonical section $\sigma\colon \R\ra \gmv$
for which $\sigma(t)$ is the Gauss-norm corresponding to the
radius~$e^t$:
\[ \abs{P(\sigma(t)}= \sup_{n\in \Z} \abs{a_n}e^{nt},\quad
\text{if $P=\sum a_n T^n$.}\]
This section~$\sigma$ is a homeomorphism onto its image $S(\gmv)$
which is called the \emph{skeleton} of~$\gmv$.

In higher dimensions, we have a similar coordinate-wise
tropicalization map
$\tau\colon \gmvg\ra\R^g$ and a section~$\sigma$ whose
image $S(\gmvg)$ is the skeleton of~$\gmvg$.

In the case of a uniformized totally degenerate Abelian variety,
one can tropicalize its universal cover and mod out by
the image of the lattice~$\Omega$. This gives a diagram:
\[
\begin{tikzpicture}[normal line/.style={->}]
\matrix (m) [matrix of math nodes, row sep=3em,
column sep=2.5em, text height=1.5ex, text depth=0.25ex]
{ \gmvg & & \R^g \\
X^\an_v & & \R^g/\Lambda \\ };
\path[->,font=\scriptsize]
(m-1-1) edge node[auto] {$ \tau $} (m-1-3);
\path[->,font=\scriptsize]
(m-1-1) edge  (m-2-1);
\path[->,font=\scriptsize]
(m-1-3) edge  (m-2-3);
\path[->,font=\scriptsize]
(m-2-1) edge node[auto] {$ \tau_X $} (m-2-3);
\end{tikzpicture}\]
where $\Lambda=\tau(\Omega)$.
Moreover, the section~$\sigma$ descends to a section~$\sigma_X$
of~$\tau_X$ whose image~$S(X^\an_v)$ is called the skeleton of~$X^\an_v$.
This is a real torus of dimension~$g$ in~$X^\an_v$ onto which~$X^\an_v$ retracts
canonically.

The proof of the following two theorems is long and difficult
and cannot be described here.

\begin{theo}[\cite{gubler2007a}, Corollary~9.9]
The canonical measure $c_1(\bar L)_v^{\dim X}$ on~$X^\an_v$
is the unique Haar measure supported by the real torus~$S(X^\an_v)$
of total mass $\deg_L(X)$.
\end{theo}

\begin{theo}[\cite{gubler2007a}, Theorem~1.3]
\label{theo.gubler-tropical}
Let $Y$ be an integral subvariety of~$X$; let $d$ be its dimension.

The image $\tau_X(Y^\an_v)$ is a union of simplices of~$\R^g/\Lambda$
of dimension~$d$.

Restricted to any of those simplices,
the direct image $(\tau_X)_*(c_1(\bar L|_Y)^{\dim Y}_v)$
on~$\R^g/\Lambda$ of the canonical measure of~$Y$
is  a positive multiple of the Lebesgue measure.
\end{theo}

\subsection{}
Given the last two theorems, \cite{gubler2007b} can complete the
proof of the Bogomolov conjecture when the given Abelian
variety has totally degenerate reduction at the place~$v$.

Indeed, in the above situation of
a generically finite map
$\phi\colon Y^m\ra W\subset X^{m-1}$, one can push
the equality of measures
$\phi_* \mu_v(Y^m)=\mu_v(W)$ to the tropicalization~$(\R^g/\Lambda)^{m-1}$.
Let $\nu_Y=(\tau_X)_*(\mu_v(Y))$, $\nu_W=(\tau_{X^{m-1}})_*(\mu_v(W))$;
these are measures on~$(\R^g/\Lambda)$ and $(\R^g/\Lambda)^{m-1}$
respectively.
Let $\psi$ be the map~$(\R^g/\Lambda)^m\ra(\R^g/\Lambda)^{m-1}$
given by $(a_1,\dots,a_m)\mapsto (a_2-a_1,\dots,a_m-a_{m-1})$.
By naturality of tropicalization, one has
$ \tau\circ \phi = \psi\circ \tau$, hence
$\psi_*(\nu_Y^m)=\nu_W$.

Let $\delta$ be a simplex of dimension~$\dim(Y)$ appearing in~$\tau_X(Y)$.
By Theorem~\ref{theo.gubler-tropical}, the restriction of
the measure~$\nu_Y$ to~$\delta$ is a positive multiple
of the Lebesgue measure. In particular, $\nu_Y(\delta)>0$.
Then $\delta^m$ is a simplex of~$\tau_{X^m}(Y^m)$ whose image by~$\psi$
is $\psi(\delta^m)$.
However, the definition of~$\psi$ shows that $\psi(\delta^m)$
has dimension~$\leq m\dim(Y)-\dim(Y)<m\dim(Y)=\dim(W)$. Indeed, $\psi$ is linear
and the simplex~$\delta$ embedded diagonally into~$\delta^m$ maps to~$0$.
By Theorem~\ref{theo.gubler-tropical},
$\nu_W$ is a sum of Lebesgue measures of $\dim(W)$-dimensional
simplices, so that $\nu_W(\psi(\delta^m))=0$.
Since $\psi_*(\nu_Y^m)=\nu_W$, it follows that $\nu_Y(\delta)=0$.
This contradiction concludes Gubler's proof of the Bogomolov conjecture
when there is a place of totally degenerate reduction.

\subsection{}
In our discussion of Manin--Mumford's conjecture over function fields,
it was necessary to take care of constant Abelian subvarieties.
They do not appear in Gubler's statement. Indeed, if an Abelian
variety has totally degenerate reduction at some place, it cannot
contain any constant Abelian subvariety. However, a general
treatment of Bogomolov's conjecture over function fields would
take them into account. A precise statement is given in the paper by~\cite{yamaki2010}, with partial generalizations  of Gubler's result
to cases where there is bad reduction, although not totally degenerate.

\bibliographystyle{mynat}
\bibliography{aclab,acl,barbados}
\end{document}